\documentclass[11pt]{article}

\usepackage[nice]{nicefrac}

\usepackage{fullpage,url,boxedminipage}
\usepackage[boxed,linesnumbered]{algorithm2e}

\usepackage{amsmath,amstext,amsfonts,amssymb}

\newcommand{\N}{\mathbb{N}}

\newcommand{\R}{\mathbb{R}}

\newcommand{\Con}{\mathcal{C}}

\renewcommand{\epsilon}{\varepsilon}
\newcommand{\maxv}[1]{{\overline{#1}}}
\newcommand{\maxeps}{\maxv{\epsilon}}

\newcommand{\erf}{\ensuremath{\mathrm{erf}}}
\newcommand{\argerf}{\ensuremath{\mathrm{argerf}}}

\newcommand{\qed}{\nobreak \ifvmode \relax \else \ifdim \lastskip<1.5em \hskip-\lastskip
\hskip1.5em plus0em minus0.5em \fi \nobreak \vrule height0.75em width0.5em depth0.25em\fi}

\begin{document}
\sloppy

  \title{Optimizing polynomials for floating-point implementation}
   \author{Florent de Dinechin \and Christoph Lauter}
\date{\today \\ ~ \\ { \small
This is LIP Research Report number RR2008-11 \\
Ceci est le Rapport de Recherches num\'ero RR2008-11 du LIP \\
Laboratoire de l'Informatique du Parall\'elisme, ENS Lyon/CNRS/INRIA/Universit\'e de Lyon \\ 46, all\'ee d'Italie, 69364 Lyon Cedex 07, France }
}

\maketitle

\begin{abstract}
  The floating-point implementation of a function on an interval often
  reduces to polynomial approximation, the polynomial being
  typically provided by Remez algorithm. However, the floating-point
  evaluation of a Remez polynomial sometimes leads to catastrophic
  cancellations. This happens when some of the polynomial coefficients
  are very small in magnitude with respects to others. In this case,
  it is better to force these coefficients to zero, which also reduces
  the operation count. This technique, classically used for odd or even
  functions, may be generalized to a much larger class of functions.
  An algorithm is presented that forces to zero the smaller
  coefficients of the initial polynomial thanks to a modified Remez
  algorithm targeting an incomplete monomial basis.  One advantage of
  this technique is that it is purely numerical, the function being
  used as a numerical black box. This algorithm is implemented
  within a larger polynomial implementation tool that is demonstrated
  on a range of examples, resulting in polynomials with less
  coefficients than those obtained the usual way.

 {\bf Keywords:} Polynomial evaluation, floating-point,
  elementary functions.
\end{abstract}

%=====================================================================
\section{Introduction}
\label{sect:intro}

Scientific and business applications need support for mathematical
functions such as $e^x$, $\sin x$, $\log x$, $\erf x$, and many
others.
Current processors offer little or no hardware for the evaluation of
such functions.  What they offer is high performance floating-point
hardware for basic operations such as addition and multiplication, or
their combination as a fused multiply-and-add.  Other mathematical
functions are usually approximated on a small domain by polynomials,
which can then be evaluated using the high-performance operators of
the processor~\cite{Muller2006,Markstein2001,CorneaHarrisonTang2002}.
This article addresses the design of such polynomial approximations.
The difficulty is to obtain floating-point implementations of
high quality with respect to both performance and precision.

\subsection{Obtaining approximation polynomials}
 
Several textbooks
\cite{Markstein2000,CorneaHarrisonTang2002,Muller2006} discuss
techniques to obtain good approximation polynomials. One may use
Taylor polynomials, Chebychev approximations, or a minimax
approximations given by Remez algorithm~\cite{Remez69,
  Muller2006}. The latter is preferred as it is theoretically the most
accurate on a domain: The minimax approximation of a function $f$ on a
domain $I$ is defined as the polynomial
$p$ %$p = \sum\limits_{i=0}^n p_i x^i$
that minimizes the infinite norm of the approximation error $\|
\nicefrac{p}{f} - 1 \|^I_\infty$.

Throughout this article we focus on
relative errors, which are more relevant to floating-point
implementations. However, a good implementations of Remez algorithm can
accommodate any kind of error through the use of an additional weight
function $w$: Remez algorithm will actually minimize $\| p \cdot w - f
\|^I_\infty$ for any $w$ provided by the user.

A few well-known tricks are used for certain classes of functions. For
instance, for an odd/even function, a straightforward application of
Remez algorithm will usually not provide an odd/even polynomial,
although an odd/even approximation polynomial has several advantages.
Firstly, it has twice as few coefficients as the minimax polynomial,
and its evaluation will therefore be more economical. Secondly, the
implementation will have the same symmetry properties as the function
implemented. Thirdly, it is more stable numerically, as the sequel
will show. The textbooks therefore advocate the use of an ad-hoc
technique (which will be reviewed in Section~\ref{sec:intro-examples})
to compute, using standard Remez algorithm, a minimax approximation
among the set of odd (resp. even) polynomials.

This article eventually offers a generalization of this ad-hoc
technique to a much larger class of functions, but it didn't start
this way.  The initial motivation of this work was to avoid a problem
that makes the evaluation of some polynomials impractical using
floating-point arithmetic.  This problem is the occurrence of
 cancellations in the evaluation. It would in particular
plague the evaluation of a straightforward minimax approximation to an
odd/even function.

\subsection{Floating-point evaluation of a polynomial}

Let us consider the Horner evaluation of a polynomial of degree $d$:
$$p(x)=a_0+x\cdot(a_1 + x\cdot(a_2+... )...) \quad .$$
It can be described as a recurrence:
$$\left\{
  \begin{array}{rcl}
    q_n & = & a_n\\
    q_{i} &=& a_{i} + x\cdot q_{i+1} \quad \text{ for } i = 0...n-1\\
    p(x) &=& q_0
  \end{array}
\right.$$

A cancellation happens when, for some values of $x$, the addition in
$a_i + x\cdot q_{i+1}$ is effectively a subtraction and the two terms
$a_i$ and $ x\cdot q_{i+1}$ are very close to each other.

Cancellations in polynomial evaluation should be avoided for two reasons.

\begin{itemize}
\item Although a cancelling subtraction is an exact operation, in the
  sense that it never involves a rounding, some digits of the result
  are no longer significant. Consider for instance, in a decimal
  format with 4-digit significand, the subtraction $1.234-1.232$. It
  is an exact subtraction that cancels three digits and returns
  $2.000\cdot 10^{-3}$: the three zeroes do not correspond to
  significant information from the initial numbers. If this result is
  then multiplied by another value --which is precisely the case in
  Horner evaluation-- this loss of information propagates to the
  result of the multiplication. As a side effect, error analysis
  becomes very difficult. This is the reason why useful theorems about the
  accuracy of Horner evaluation \cite{BolDau04b} have hypotheses that
  ensure that cancellation will not occur.

\item Fast double-double \cite{Knuth97} and
  triple-double~\cite{Lauter2005LIP:tripledouble} operations can only  be
  used if the absence of cancellation can be proven beforehand.
  Otherwise, much slower versions have to be used~\cite{Finot-thesis}.
\end{itemize}

The second reason leads us to avoid any cancellation.  Usually this
means that $x\cdot q_{i+1}$ should remain smaller than
${a_i}/{2}$ or that both always have the same sign -- see
Section~\ref{sec:complete-alg} for the actual criterion used.

Other evaluation schemes exist and are better suited to current
superscalar machines, from the family of Estrin
schemes~\cite{CorneaHarrisonTang2002,Muller2006,DinLauMul2006:log}
to Knuth/Eve transformations \cite{Knuth97,Revy2006}.  The issue of
cancellation is relatively independent of the choice of evaluation
scheme, since it is mostly related to the orders of magnitude of the
coefficients. More specifically, these other schemes can also be
expressed as a tree of operations involving additions and partial
polynomials, and  the techniques developed for Horner
should adapt to other schemes as well.

\subsection{Contributions}

It is quite simple to detect the possibility of cancellations in
a polynomial evaluation scheme.
We present an algorithm that starts with a standard minimax polynomial
for a function, searches its Horner evaluation for cancellations, then
tries to remove these cancellations by setting to zero the offending
coefficients: if $a_{i} + x\cdot q_{i+1}$ may cancel, then set $a_i$
to zero. The algorithm then iterates using a modified Remez algorithm
that forces some coefficients to zero.

It turns out that this algorithm numerically discovers the zero
coefficients of odd/even functions, but also in many more cases when
textbook recipes do not apply -- examples will be given in
Section~\ref{sec:intro-examples}. Even for odd/even functions, it
widens the implementation space, as it explores implementations that
are not strictly odd/even. It thus provides practical answers to
questions such as: Is this odd polynomial  really the optimal one
with respect to my optimality criterion?

Moreover, this algorithm is purely numerical: all it requires is a
black-box implementation of the function and its first two
derivatives. This implementation has to return an enclosure of these
functions in multiple precision, up to an accuracy sufficient with
respect to the target precision. It may therefore be applied to
arbitrary compound functions, but also to functions defined by
integrals or iterative processes -- in such cases the enclosures of
the first two derivatives $f$ may have to be computed numerically out
of those of $f$. We illustrate this in
Section~\ref{sec:intro-examples} on the example of $\argerf =
\erf^{-1}$, which we compute by a Newton-Raphson-iteration out of on
the function $\erf$.

This algorithm may even be used for a class of applications studied by
Remez~\cite{Remez69}, then Dunham~\cite{Dunham}, where a function is
defined by a large set of points obtained from some physical
experiment. A linear interpolation between these points allows one to
recover the function, but one may want a more compact polynomial
implementation. Provided the set of points is large enough and
accurate enough, the technique presented here should apply to such
problems, too.

Our algorithm has been implemented using the Sollya
tool\footnote{\url{http://sollya.gforge.inria.fr/}}, in the Sollya
scripting language.  Sollya provides all the necessary building
blocks.  Internally, it uses multiple precision interval arithmetic
based on the MPFR\footnote{\url{http://www.mpfr.org/}} library. Based
on this, Sollya evaluates composite functions with faithful rounding
to any precision. Functions written outside Sollya can also be linked
dynamically. In addition, Sollya also provides a fast but nevertheless
high-quality implementation of the infinite norm, and also a slower,
certified one~\cite{LauterChevillard2007}. Finally, Sollya integrates
a Remez algorithm (presented in Section~\ref{sec:remez}) which is more
flexible than other implementations, for instance the minimax
algorithm in the Maple \texttt{numapprox} package.

% Besides, the validity of the resulting polynomial can be
% verified a posteriori with safe algorithms. In our case it suffices to
% certify a bound on the approximation error of the
% polynomial~\cite{LauterChevillard2007}.

\subsection{A complete implementation chain}

The algorithm presented in this article  provides a polynomial with real coefficients. What one
eventually needs for a floating-point implementation is a polynomial
with machine-representable coefficients: In the applications we
target~\cite{DinLauMul2006:log}, coefficients should be double-precision numbers,
double-doubles, or triple-doubles
\cite{Bailey2001,Lauter2005LIP:tripledouble}. Techniques are known to
obtain a minimax polynomial among the set of polynomial with
such machine-representable coefficients~\cite{BrisebarreChevillard2007}.
These techniques have to be initialized with a good polynomial: the
polynomial obtained by our algorithm fills this role.  This final
transformation of our polynomial is therefore out of the scope of the
present article, although some of the examples provided in
Section~\ref{sec:examples} have gone through this final step, which is
also automatic. The appropriate precision for each addition and
multiplication step in the Horner scheme is chosen automatically.  Let
us just mention that one needs to ensure that this final step keeps
our zero coefficient, which is easy within the framework
of~\cite{BrisebarreChevillard2007}. One also has to ensure that this
final step does not introduce a possibility of cancellation. This can
be checked a posteriori, and in practice the check always succeeds,
because the algorithm of \cite{BrisebarreChevillard2007} preserves
the floating-point binade of the coefficients.

Finally one has to evaluate or certify a bound on the accumulated
rounding error when this final polynomial is evaluated. This can be
done using the Gappa tool~\cite{DinLauMelq2005:gappa}, and is also
automated, but out of scope of this article.

To sum up, the algorithm presented in this article is integrated in a larger
Sollya program that is capable of automatically generating certified C code for
Horner's scheme using expansion arithmetic from double to
triple-double precision~\cite{DinLauMul2006:log}. Let us now focus on this algorithm.

\section{Minimax approximation on an incomplete monomial basis}
\label{sec:remez}

This section first revisits the current textbook methodology for
approximating a function with a polynomial, then introduces several
useful generalizations and formalizations.

\subsection{Introductory examples}
\label{sec:intro-examples}

Let us consider the example of an implementation of the sine function
between $-\nicefrac{\pi}{64}$ and $\nicefrac{\pi}{64}$ (this interval is typical after a
table-based argument reduction \cite{Markstein2000,Muller2006}). Let
us assume that a target accuracy of $2^{-60}$ is wanted, as would
be the case for a quality double-precision implementation.

A first option is to use a Taylor approximation. As the sine
function is odd, the corresponding series has only odd monomials.
With degree 7, the relative approximation error is roughly $2^{-53}$. With
degree 9, the error is smaller than $2^{-68}$. This illustrates how
discrete the quality of the approximation usually is with respect to
the degree.
%  p=taylor(sin(x), 7,0);
% e=dirtyinfnorm(p/sin(x)-1, [0,Pi/128]); log2(e);
% 
% A second option is to use a Chebyshev approximation. TODO regarder le Muller

Another option is to use a minimax polynomial, provided by Remez algorithm.
This will provide us with a polynomial that is not odd: its coefficients are 
% p = remez(1,[|1,...,7|],[0,Pi/64],1/sin(x));
\[
\small
\begin{array}{lll}
a_1= 1.00000000000000000004553862129419953814366183346717373\\
a_2= -0.874967378163014390017316615896238907007870683208826576 \cdot 10^{-16}\\
a_3= -0.166666666666639297309612035148824100626593835839529139\\
a_4= -0.317973607302440662105040928632951877918380396824236704 \cdot 10^{-11}\\
a_5= 0.833333350548021113401528637698582467442968048220758103 \cdot 10^{-2}\\
a_6= -0.454284495616307678859183047154098346831822774513536551 \cdot 10^{-8}\\
a_7= -0.198362485524232245861352857470565050846399633158981546 \cdot 10^{-3}
\end{array}
\]

The approximation error is now smaller than $2^{-64}$, and this polynomial
fits our accuracy target. However, this full
polynomial of degree 7 will require 7 multiplications and 6 additions
when evaluated by Horner, whereas the Taylor of degree 9, considered
as a degree-4 polynomial in $x^2$ multiplied by $x$, requires only 6 multiplications and
4 additions, and provides a much better approximation error
($2^{-68}$). 

This is a good enough reason to prefer the degree-9 Taylor, but there
is another critical problem with the minimax polynomial: when
evaluated in finite precision, the alternation of comparatively very
small and very large coefficients  will lead to large cancellations.

Still, the approximation error of the Taylor polynomial is not
balanced on the interval. One may think of zeroing out the smaller
(even) coefficient of the minimax, but this provides a much worse approximation
% p=subpoly(p,[|1,3,5,7|]);
error ($2^{-49}$). Textbooks by Muller and Markstein therefore suggest
the following better recipe for finding an odd polynomial with a
better balanced error: perform the change of variable $X=x^2$, and
compute a degree-4 minimax on the transformed function
$\sin(\sqrt{X})/\sqrt{X}$, with a suitable weight function. This will
provide a degree-4 polynomial in $X$. Multiplying the corresponding
degree-8 polynomial in $x$ with $x$ yields an odd polynomial that
approximates $\sin(x)$ well: the error is smaller than $2^{-60}$. This
polynomial evaluates in only 6 multiplications, and without
cancellation issues -- this last point as a side effect, not mentioned
in textbooks. Besides, a good choice of the weight function  ensures that this
polynomial is the minimax among the set of odd polynomials of
degree 9.

This recipe is more or less the state of the art for finding an
approximation polynomial, and the improvements over it have mostly
consisted in obtaining minimax approximations over the set of
polynomial with machine-representable coefficients
\cite{BrisebarreChevillard2007}. Let us now try to apply it to more
complex functions.

\begin{itemize}
\item $\cos(x^2)$ has a Taylor approximation of the form
  $1-1/2x^4+1/24x^8$ with only coefficients for powers of $x^4$. The
  recipe above can be accommodated to such cases.

\item $\sin(x) + \cos(x^2)$ has the following Taylor approximation:\\
  $1 + x - 1/6 x^3 - 1/2 x^4 + 1/120 x^5 - 1/5040 x^7 + 1/24 x^8 +
  1/362880 x^9...$.  Here the zero coefficients -- which will probably
  lead to cancelling coefficients in minimax approximation -- are the
  coefficients of degree $2+4k$ for $k$ integer. This structure can be
  predicted fairly straightforwardly from the sums of the Taylor series
  of $\sin(x)$ and $\cos(x^2)$. However, the textbook recipe does not
  extend straightforwardly to this function: providing a minimax
  polynomial with this coefficient structure is not obvious. One may
  apply the recipe to $\sin(x)$, then to $\cos(x^2)$, then add the
  resulting polynomials coefficient-wise. This will provide a
  polynomial with the wanted coefficient structure, but there is no
  reason why it should be the one with minimal error.

\item $e^{\sin(x) + \cos(x^2)}$ has a Taylor series with a single zero
  coefficient at degree 3. Again, this can be predicted analytically
  or using computer algebra, but it is getting less and less
  intuitive. And again, can the reader provide a change of variable
  allowing to compute a minimax polynomial for this function among the
  set of polynomial whose degree-3 coefficient is zero?

\end{itemize}

Looking back at the Taylor series for $\sin(x) + \cos(x^2)$, we also
remark that even its non-zero coefficients will present risks of
cancellation as the degree augments. Any methodology leading to an
implementation has to study these cancellation issues at some
point. Our approach is to attack the cancellation issues first by
imposing the value 0 to some coefficients. It happens that such a
methodology will, as a side effect, zero out the coefficients which
are zero in the Taylor series, and thus minimizes the number of
operations required for the evaluation, just as the textbook recipe, but for more functions. 

The core of this technique is a modified Remez algorithm that is
able to compute a minimax polynomial on an incomplete monomial basis,
which we present now.

% Remarque que Muller optimise le nombre d'operations et au passage obtient un polynome evaluable, nous, on veut evaluer et l'assurer et 
% on obtient l'optimisation du nb d'operation au passage

% Attaquons des fonction plus compliquees + remarque sur Taylor => evaluable relié au nb d'operation minimal

% exp(cos(x^2) + 1)

% saut conceptuel, exp(sin x - cos x^2) 

% => ici découvrir que les recettes ne suffisent plus (typiquement sin x - cos x^2) => Remez à trous

\subsection{Modifications to Remez algorithm}

Refer to~\cite{Muller2006} for a good introduction to the Remez
algorithm~\cite{Remez69, Cheney66}. Modifying this algorithm for an
incomplete monomial basis $\left \lbrace x^{i_0} , x^{i_2}, \dots,
  x^{i_n} \right \rbrace$ is fairly straightforward. The core of Remez
algorithm solves a linear interpolation problem on $n + 2$ points.  In
the general case, the matrix of this linear system is a Vandermonde
matrix. For an incomplete basis, the matrix simply only contains columns
corresponding to the given monomials. This
possibility is already present in the initial formulation of the problem by
Remez~\cite{Remez69, Cheney66}.

% Typically the system matrix $M$ becomes the
% matrix of the $\left( x^{i_k}_j \right)_{kj}$. 

\subsection{Theoretical issues}

Remez algorithm works iteratively. In order to ensure that the
polynomial is really improved in each step, classical Remez algorithm
requires the so-called Haar condition~\cite{Remez69, Cheney66}: the
determinant of the formal matrix of the linear system does not vanish
for any choice of approximation points in the interval $I$. This Haar
condition is always satisfied for Vandermonde matrices~\cite{Remez69},
hence for approximations on the complete monomial basis.

Unfortunately, the Haar condition is not fulfilled in general for
incomplete bases, in particular if the interval comprises
zero~\cite{Remez69, Cheney66}. Therefore the convergence of our
modified Remez is not proven. Furthermore, when it converges, there is
currently no proof that it converges to the expected minimax
polynomial.

These issues are under investigation. To the best of our knowledge,
there exists some literature addressing them, sometimes dating back
to Remez himself~\cite{Remez69,Schaback70}. However there is no recent
work on a multiple-precision implementation of a minimax without Haar
condition.

\subsection{The modified Remez algorithm in practice}
\label{sec:remez-pract}
In practice our modified Remez algorithm, implemented in the Sollya
tool, satisfies our needs. The algorithm typically runs even without
Haar condition until it converges to some polynomial~\cite{Remez69},
or until it can no longer exhibit $n + 2$ appropriate interpolation
points. When it works, the returned polynomial might not be the actual
best approximation polynomial. Nevertheless, it is possible to verify
a posteriori whether or not this polynomial satisfies the target
approximation error bound. This is enough to validate the
implementation of a function. Experimentally, comparisons with minimax
on the complete basis, as well as observed discrete jumps of errors
with respect to the degree, suggest that the polynomial is not far
from the optimal. Actually, the optimality can be proven a posteriori
using a modification of the polytope exploration described
in~\cite{BrisMulTiss2006}, but this is very costly and out of scope of this paper.

In addition, we use the modified Remez algorithm with care. For
example, it is possible to ensure the Haar condition in cases when the incomplete
monomial basis is even/odd: If the
approximation interval $I$ comprises zero, just take only the positive
(or negative) half of the interval.

\section{Computing cancellation-free approximation polynomials}
\label{sec:complete-alg}

The complete algorithm implementing our approach is sketched in the
listing below.  
 It first tests if the approximation
 polynomial in the complete basis is cancellation-free. If the test
 fails, all monomials on which cancellations occur are removed from the
 basis. An approximation polynomial is then computed in the resulting,
 incomplete basis. The algorithm iterates on increasing degrees until
 the target error bound is fulfilled or an iteration limit is reached.

\begin{algorithm}[H]\label{algo}
\SetLine
\SetCommentSty{textsl}

\SetKwData{Cancellationfree}{cancellationfree}
\SetKwFunction{Guessdegree}{guessdegree}
\SetKwFunction{Remez}{remez}

\KwIn{a function $f \in \Con^2$, a domain $I$, a target error $\maxeps \in \R^+$, an iteration limit $l \in \N$}
\KwOut{a polynomial $p$, cancellation-free in Horner's scheme, such that $\| \nicefrac{p}{f} - 1 \|^I_\infty \leq \maxeps\quad$ $\bot$ if no such polynomial can be found}
\BlankLine

$n \leftarrow$ \Guessdegree{$f$,$I$,$\maxeps$}$-1$ \tcc*{Compute initial guess of the required degree}
\Repeat(\tcc*[f]{Iterate until required degree $n$ is found}){$\epsilon \leq \maxeps$}{
$n \leftarrow n + 1$\;
$p^* \leftarrow$ \Remez{$f$, $I$, $\left \lbrace x^0, \dots, x^n \right \rbrace$} \tcc*{Compute polynomial in the complete basis} \nllabel{remez1}
$\epsilon \leftarrow \| \nicefrac{p^*}{f} - 1 \|^I_\infty$ \tcc*{Compute bound on approximation error}
}
$k \leftarrow 1$\tcc*{Iteration count}
\Repeat(\tcc*[f]{Iterate until appropriate cancellation-free polynomial found}){$k > l$}{
\Cancellationfree $\leftarrow$ true\;
$B \leftarrow \left \lbrace x^n \right \rbrace$ \tcc*{$B$ will hold basis of non-cancelling additions}
$q_n \leftarrow p^*_n$\;
\For(\tcc*[f]{Statically simulate steps of Horner's scheme}){$i \leftarrow n-1$ \emph{\KwTo}$0$}{
$\overline{\alpha} \leftarrow \sup \left \lbrace x \cdot q_{i+1}(x) \vert x \in I \right \rbrace$;
$\quad \underline{\alpha} \leftarrow \inf \left \lbrace x \cdot q_{i+1}(x) \vert x \in I \right \rbrace$;
$\quad \alpha \leftarrow \max \left( \left \vert \underline{\alpha} \right \vert , \left \vert \overline{\alpha} \right \vert \right)$\;
\eIf(\tcc*[f]{No cancellation?}){$\left( \alpha \leq \nicefrac{1}{2} \cdot \left \vert p^*_i \right \vert \right)$ or 
$\left( \underline{\alpha}, \overline{\alpha}, p^*_i \right.$ have the same sign$\left. \right)$}{
$B \leftarrow B \cup \left \lbrace x^i \right \rbrace$ \tcc*{No cancellation, add monomial to basis}
}{
\Cancellationfree $\leftarrow$ false \tcc*{Cancellation occurs}
}
$q_{i} \leftarrow p^*_i + x \cdot q_{i+1}$ \tcc*{Statically simulated step in Horner's scheme}  \nllabel{symbolic}
}
\lIf{\Cancellationfree}{\Return $p^*$}\tcc*{Polynomial in full basis is cancellation-free}\Else{
$p \leftarrow$ \Remez{$f$, $I$, $B$}\tcc*{Compute approx. polynomial in incomplete basis $B$} \nllabel{remez2}
$\epsilon \leftarrow \| \nicefrac{p}{f} - 1 \|^I_\infty$\tcc*{Compute associated error bound}
\lIf(\tcc*[f]{Error of approx. polynomial in incomplete basis okay?}){$\epsilon \leq \maxeps$}{\Return $p$}
\Else{
$n \leftarrow n + 1$\tcc*{Increase the degree $n$ by $1$}
$p^* \leftarrow$ \Remez{$f$, $I$, $\left \lbrace x^0, \dots, x^n \right \rbrace$}\tcc*{Compute polynomial in complete basis}
$\epsilon \leftarrow \| \nicefrac{p^*}{f} - 1 \|^I_\infty$\tcc*{Compute associated error bound}
} 
}
$k \leftarrow k + 1$\;
}
\Return $\bot$\tcc*{No polynomial has been found in $l$ iterations}
\caption{The complete approximation algorithm}
\end{algorithm}

The \texttt{guessdegree}  Sollya function  returns a guess of the degree needed for
approximating a function $f$ in a domain $I$ with error less than
$\maxeps$. It is based on the first
iteration of the Remez algorithm combined with a bisection search.

The assignment $q_{i} \leftarrow p^*_i + x \cdot q_{i+1}$ at line
\ref{symbolic} of the algorithm must be understood as a symbolic
multiplication of the (symbolic) polynomial $q_{i+1}$ by a free variable
$x$ and a symbolic addition of $p^*_i$. 

As mentioned in section \ref{sec:remez-pract}, care is needed in the
calls to the Remez algorithm. Lines \ref{remez1} and \ref{remez2} of
the algorithm above refer to sub-routines that represent a few
hundreds of lines of code in the Sollya scripting language.

All computations of infinite norms, infima and suprema occurring
in the algorithm can be performed by uncertified, quick-and-dirty
numerical algorithms. This permits to use a black-box implementation
for the function $f$. Certification of the infinite norm of the final
polynomial  can be performed a
posteriori~\cite{LauterChevillard2007}.

An iteration limit $l$ is taken in input of the algorithm. We introduce
it for ensuring the termination of the algorithm on degenerated
problems that we believe to exist. Currently we can neither prove the
convergence of the algorithm nor exhibit an example on which it would
not terminate. In mathematical terms, we are currently unable to
analyze whether the set of cancellation-free polynomials (as defined
above) of arbitrary degree is dense for any function and approximation
interval. In practice, the iteration limit is not reached. 

% Code de moulinette

% Remarque boite noire 

\section{Examples}
\label{sec:examples}
\subsection{From a function to a C implementation}
Let us first reconsider the function $f$ defined by $f(x)
= e^{\sin x - \cos x^2}$ in the domain $I = \left[ -2^{-8}; 2^{-8}
\right]$. The function is to be approximated by a polynomial with a
relative error less than $2^{-90}$. Our algorithm chooses the monomial
basis $\left \lbrace x^0, x^1, x^2, x^4, x^5, x^6, x^7, x^8, x^9 \right
\rbrace$, leaving out the monomial $x^3$. The implemented polynomial
$p(x) = \sum\limits_{i=0}^0 p_i x^i$ is given by
\begin{eqnarray*}
p_0 & = & 119383704169626743428469396878343 \cdot 2^{-108} \\
p_1 & = & 29845926042406685857117349204375 \cdot 2^{-106} \\
p_2 & = & 119383704169626743428436621385363 \cdot 2^{-109} \\
p_4 & = & 4970345142530923 \cdot 2^{-55} \\
p_5 & = & 358969371405011 \cdot 2^{-51} \\
p_6 & = & 6516674741954513 \cdot 2^{-56} \\
p_7 & = & 589077943038783 \cdot 2^{-57} \\
p_8 & = & 5559725200690211 \cdot 2^{-59} \\
p_9 & = & 5320394595779079 \cdot 2^{-58} 
\end{eqnarray*}
As can be easily verified by computing $\| \nicefrac{p}{f} - 1
\|^I_\infty$, its error is bounded by $2^{-90.4}$. It can be evaluated
using Horner's scheme with IEEE 754 double and double-double
arithmetic with an round-off error of less than $2^{-93.6}$. This
bound has been shown with the Gappa tool~\cite{DinLauMelq2005:gappa},
and in the process Gappa formally shows a posteriori that no
cancellation may occur.

\subsection{Adaptation of the monomial basis  to the target accuracy}
Let us consider now the function $f$ defined by $f(x) = e^{\cos x^2 +
  1}$ in the domain $I = \left[ -2^{-8}; 2^{-5} \right]$. The table
below shows the monomial basis chosen by our algorithm, as a function of
the  relative approximation accuracy target.

\begin{center}
   \renewcommand{\arraystretch}{1.2}
\begin{tabular}{|l|l|}
\hline
accuracy target & monomial basis \\
\hline \hline 
$2^{-40}$ & $\left \lbrace x^0, x^4 \right \rbrace$ \\
\hline
$2^{-50}$ & $\left \lbrace x^0, x^4, x^8 \right \rbrace$ \\
\hline
$2^{-60}$ & $\left \lbrace x^0, x^4, x^8 \right \rbrace$ \\
\hline
$2^{-70}$ & $\left \lbrace x^0, x^4, x^8, x^{12} \right \rbrace$ \\
\hline
$2^{-80}$ & $\left \lbrace x^0, x^4, x^8, x^{12} \right \rbrace$ \\
\hline
$2^{-90}$ & $\left \lbrace x^0, x^4, x^8, x^{12} \right \rbrace$ \\
\hline
$2^{-100}$ & $\left \lbrace x^0, x^4, x^8, x^{12}, x^{13}, x^{14}, x^{15} \right \rbrace$ \\
\hline
$2^{-110}$ & $\left \lbrace x^0, x^4, x^8, x^{12}, x^{16} \right \rbrace$ \\
\hline
$2^{-120}$ & $\left \lbrace x^0, x^4, x^8, x^{12}, x^{16}, x^{17}, x^{18} \right \rbrace$ \\
\hline
\end{tabular}
\end{center}

\subsection{An example of black-box function}
The following example illustrates that the algorithm works well on
black-box functions. The function $\argerf = \erf^{-1}$ is not
currently available in the MPFR
library, and the Sollya tool is
based on MPFR. It is nevertheless relatively easy -- although probably very inefficient -- to implement
using a Newton-Raphson-iteration on the function $\erf$. Its first two
derivatives can be computed based on the code for $\argerf$ itself.
We have implemented such a black-box function $\argerf$ and have
dynamically bound it to Sollya.  

For an approximation polynomial in the domain $I = \left[
  -\nicefrac{1}{4}; \nicefrac{1}{4} \right]$ and a relative error
bound of $2^{-60}$, our algorithm finds the monomial basis $\left
  \lbrace x^1, x^3, x^5, x^7, x^9, x^{11}, x^{13}, x^{15}, x^{17},
  x^{19} \right \rbrace$. As $\argerf$ is given as a binary
executable, there is no formal knowledge that it is an odd function.
Nevertheless, the algorithm chooses an odd monomial basis for this odd
function. The polynomial $p(x) = \sum\limits_{i=0}^{19} p_i x^i$
implemented using IEEE 754 doubles and double-doubles is given by
\begin{eqnarray*}
p_1 & = & 71899270015270848535577833907197 \cdot 2^{-106} \\
p_3 & = & 37646369746407330411070885976913 \cdot 2^{-107} \\
p_5 & = & 2297847774298601 \cdot 2^{-54} \\
p_7 & = & 3118369096730189 \cdot 2^{-55} \\
p_9 & = & 2340416807028733 \cdot 2^{-55} \\
p_{11} & = & 7455281238343373 \cdot 2^{-57} \\
p_{13} & = & 3086390951797773 \cdot 2^{-56} \\
p_{15} & = & 5269462590206135 \cdot 2^{-57} \\
p_{17} & = & 8758767795225423 \cdot 2^{-58} \\
p_{19} & = & 5369190506948897 \cdot 2^{-57} 
\end{eqnarray*}
Its approximation error is bounded by $2^{-62.9}$. The polynomial can
be evaluated using Horner's scheme and double and double-double
arithmetic. The corresponding evaluation error is bounded by
$2^{-62.4}$.

\section{Conclusion and future works} \label{conclusion}

Obtaining a floating-point C implementation of an arbitrary function
has never been so automatic. The tools presented in this article are
able to provide a high quality implementation of reasonably complex
functions functions in a few seconds. Some of the a-posteriori
validation steps may require many hours, though.

In the process of developing these tools, much insight has been gained
in the issue of cancellation in polynomial evaluation, and its
relationship to the monomial basis in which polynomials should be
searched. The proposed modifications to Remez algorithm are simple in
principle but complex in the details, and in the supporting
theory. More solid mathematical foundations remain to be built, or
rediscovered from ancient literature -- the reference books on the
subject date back to the '60s. Meanwhile, we have demonstrated an
implementation that works in practice. 

Of course, expertise remains necessary in the implementation of an
elementary function. A clever argument reduction can bring in accuracy
and performance improvements by orders of magnitude above polynomial
approximation alone, and this is totally function-dependent. Even in
this case, the presented tool will be very helpful in the hands of  the
expert. Modern table-based argument reduction techniques are often
parameterized, and lead to wide trade-offs. Choosing the best set of
parameters may be very tedious, and performance-wise, it will be very
dependent on the target machine. The tools presented here will help
the designer navigate these trade-offs quickly. They have been used (at
various stages of maturity) to build some of the functions in the
CRLibm library\footnote{\url{http://lipforge.ens-lyon.fr/www/crlibm/}}.

Another, longer-term application of this kind of tool is the
construction at compile-time of ad-hoc polynomial evaluators for
composite functions appearing in application code. Many problems
remain to be solved in this context, in particular the determination
of a pertinent input interval -- this is a difficult problem of static
analysis in compilation. This approach could lead to higher
performance, but also better accuracy than the composition of
elementary functions currently used.

% Remez à trous, issues etc. modeste
% Convergence ?

\bibliographystyle{plain} 
%\bibliography{arenaire.bib,christoph.bib,arith.bib,elem-fun.bib,f2d.bib,guillaume.bib}

\end{document}